# Monitoring & Mitigation of Delayed Voltage Recovery using µPMU Measurements with Reduced Distribution System Model

Amarsagar Reddy Ramapuram Matavalam, *Student Member, IEEE*, Ramakrishna Venkatraman, *Student Member, IEEE* and Venkataramana Ajjarapu, *Fellow, IEEE.*

***Abstract*—This paper proposes a new method to monitor and mitigate fault induced delayed voltage recovery (FIDVR) phenomenon in distribution systems using µPMU measurements in conjunction with a Reduced Distribution System Model (RDSM). The recovery time estimated from a dynamic analysis of the FIDVR is used to monitor its behavior and a linear optimization is formulated to control air conditioner loads and DER reactive power injection to mitigate the FIDVR severity. The RDSM is made up of several sub-models, each of which is analogous to the Composite Load Model (CLM) with selected parameters. The linear formulation in combination with the RDSM reduces the computation time, enabling online execution. Simulated µPMU measurements from the IEEE 37 node distribution system connected to the IEEE 9 bus system under various fault scenarios are used to evaluate the proposed methodology. The resulting mitigation schemes are validated using combined transmission-distribution system simulations, thereby demonstrating that µPMU measurements along with the RDSM enable FIDVR mitigation by optimal control of reactive power injection from DERs with minimal load disconnection.**

## I. Introduction

IN today's ever evolving power grid, advanced monitoring and control schemes to mitigate abnormal grid behavior such as short term voltage instability are vital for the reliability. The phenomenon of short term voltage stability deals with the behavior of the power system in the first few seconds after a disturbance. A special case of interest is the Fault Induced Delayed Voltage Recovery (FIDVR) phenomenon which occurs in regions where the 1ϕ induction motor (IM) load portion is more than 30% [1]-[2]. FIDVR is a precursor to short term voltage instability since the generator excitation and the transmission lines are stressed due to motor stalling, thus increasing the risk of cascading. FIDVR is mainly observed in systems with a moderate proportion of 1ϕ IM loads, which are present mainly in air conditioner (A/C) loads. After a large disturbance (fault, etc.), these motors can stall and draw ~6 times their nominal current, leading to the depression of the system voltage for several seconds (>15 sec).

Two types of methodologies have been proposed in literature to mitigate the FIDVR phenomenon – supply side methods (injection of dynamic VARs via SVC, etc.) and demand side methods (disconnection of loads using measurements, etc.). Utilities usually employ the supply side solution by determining the amount and location of the SVCs and STATCOMs during the offline planning phase [3][4]. These methods use contingency sets along with extensive time domain simulations to ensure that FIDVR is mitigated over a wide range of operating conditions. The widespread adoption of Phasor Measurement Units (PMUs) by utilities has led to the development of measurement based methods to estimate the severity of FIDVR in real-time and take appropriate control actions to prevent further voltage reduction [5][6][7].

Until recently, distribution systems (DS) have lacked high-quality real-time measurement data. There has been a compelling motivation for using advanced measurement data from accurate, high resolution devices in distribution networks [8]. High-precision micro phasor measurement units (µPMUs), when tailored to the particular requirements of power distribution, can support a range of monitoring, diagnostic and control applications [8]. They can also enable a new approach for managing distribution systems, particularly in the presence of significant renewable penetration [9] and can revel phenomenon that were not usually thought to occur in distribution systems. In fact, it was recently shown from µPMU measurements that FIDVR occurred more frequently in distribution systems than transmission systems (TS) [10].

To mitigate FIDVR in distribution systems, [11][12][13] have proposed utilizing the reactive support from DER inverters based on voltage reduction at the inverter. However, as FIDVR phenomenon is driven by the load dynamics, targeted load control in regions with large motor stalling will lead to a faster recovery. This approach is adopted in this paper where we demonstrate that the µPMU measurements provide sufficient visibility to identify and localize motor stalling in distribution systems. Furthermore, by analyzing the dynamics of FIDVR, we are able to estimate the recovery time from measurements to enable improved mitigation schemes by targeted control of A/C smart thermostats and DERs. These targeted schemes are shown to mitigate FIDVR with lesser load control than widespread disconnections throughout the system.

A. R. Ramapuram Matavalam, R. Venkatraman and V. Ajjarapu are with the Department of Electrical and Computer Engineering, Iowa State University, Ames, IA 50011 USA. (e-mail: amar@iastate.edu , rvenkat@iastate.edu & vajjarap@iastate.edu).

The authors were supported from National Science Foundation grants and Department of Energy grants and are grateful for their support.

## II. Analysis & Recovery Time Estimation of FIDVR

In order to study the FIDVR phenomenon, the Composite Load Model (CLM), which is one of the most comprehensive dynamic load models, has been developed by Western Electricity Coordinating Council (WECC) [14]. This model aggregates the various loads in a region into static loads, 3ϕ IM (also referred as motor-A, B, C), and 1ϕ IM (also referred as motor-D), representing the residential A/C loads. The overall structure of the composite load model is shown in Fig. 1. The 1ϕ IM model represents the A/C compressor motor, thermal relay, and contactors and is the main reason why the FIDVR occurs. Depending on the load voltage, the 1ϕ IM operates either in 'running' or 'stalled' state.

The 1ϕ IM is in the running state for normal operating voltage and when the voltage goes below the stall voltage for a time greater than the stall time, the 1ϕ IM goes into the stalled state. In the stalled state, the active power demand is ~3 times the nominal amount and the reactive power demand is ~6 times the nominal amount [4][14] compared to the normal 'running' state. This large increase in the reactive power demand is the reason why the voltage at the load drops during stalling. This power demand is naturally reduced via thermal protection and takes around 10-15 seconds to operate. Despite the recovery, the concern is that the sustained low voltages can lead to events such as generator exciters reaching limits or disconnection of DG inverters that can initiate cascading phenomenon [15].

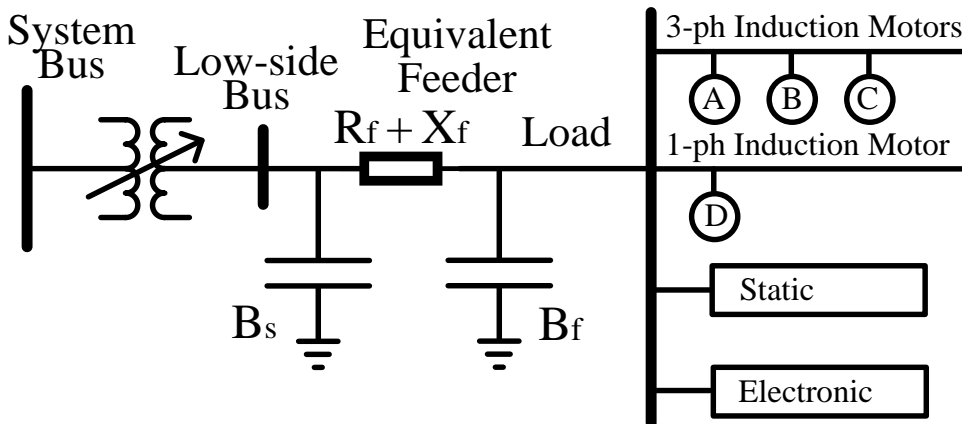

Fig. 1. Structure of the Composite Load Model [14]

Fig. 2 shows the simplified structure of the composite load model during stalling of the 1ϕ motor with the equivalent feeder admittance denoted by $\mathbf{Z_f}$. The μPMU is present at the node before $\mathbf{Z_f}$ and measures the voltage & load current which can be used to calculate the load voltage $\boldsymbol{V_L}$. As the thermal relay dynamics is much slower compared to the dynamics of the 3$\phi$ IM, the dynamics of the 3$\phi$ IM can be neglected for this particular phenomenon. The 3ϕ motor, electronic loads, static loads and DER are represented by admittance $\mathbf{Y_{3ESDER}}$ and the stalled 1ϕ motor is represented by admittance $\mathbf{Y_{stall}}$. The admittance $\mathbf{Y_{3ESDER}}$ is a function of $\boldsymbol{V_L}$ in order to account for the dynamics of the 3ϕ motor and is not constant with time. After a severe fault, the stalled 1ϕ IM admittance is given by $\mathbf{Y_{stall}} = \mathrm{G_{stall}} - \mathrm{j} \cdot \mathrm{B_{stall}}$. The fraction of 1ϕ IM connected after stalling is determined by the fraction $\mathrm{f_{th}}$ which is the output of the thermal relay. The thermal relay block diagram is shown in Fig. 3, where the thermal power dissipated in the motor given by $\mathrm{P_{th}}$ (equal to $V_L^2 \cdot G_{stall}$), $T_{th}$ is the thermal relay time constant and $\theta$ is the motor temperature estimated by the relay [14].

To validate that the load admittance can indeed capture the load behavior during FIDVR better than voltages, Fig. 4 plots the voltages and the load susceptance for moderate and severe FIDVR. The voltage waveforms have oscillations due to generator and other system dynamics. In comparison, the load susceptance is nearly flat as the generator dynamics have a much smaller impact on the load susceptance. The susceptance of the delayed voltage scenario has a sudden rise due to the stalling of the 1$\phi$ IMs. This sudden rise can be used as a reliable indicator of the FIDVR phenomenon [16].

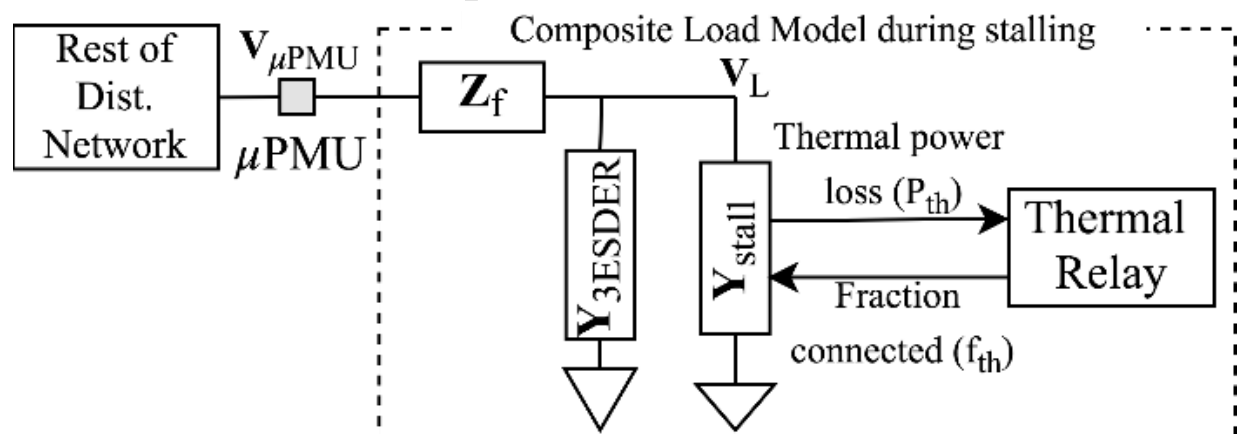

Fig. 2. Simplified composite load model during FIDVR

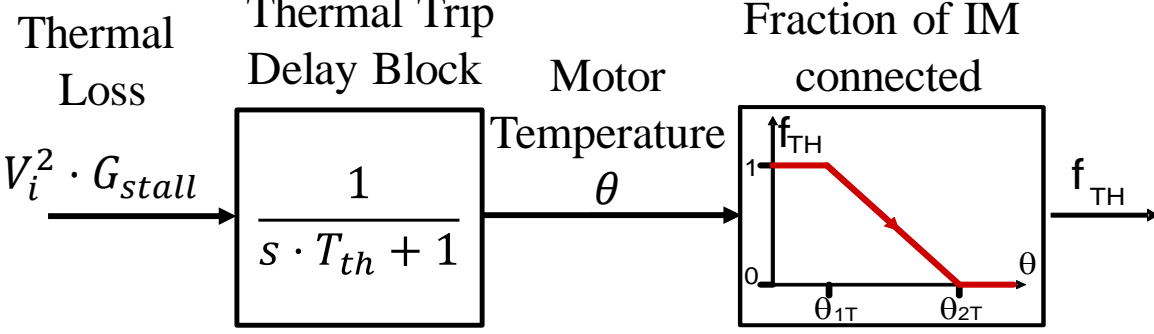

Fig. 3. The thermal relay dynamics of the 1ϕ IM [14]

Additionally, the susceptance for the delayed voltage scenario can be split into two parts – a flat region and a monotonically decreasing region. The flat region corresponds to the time to initiate the thermal tripping of 1ϕ IM ($t_1$) and the region where the susceptance reduces which corresponds to the time taken to complete the thermal tripping of 1ϕ IM ($t_2$). It is much easier to distinguish between these phases of operation from the susceptance plots compared to the voltage plots as the oscillations and other phenomenon can mask the exact time of transition. These observations and modelling assumptions lead to the admittance based representation of the composite load model [16]. Analysis of this simplified model along with the thermal relay dynamics is discussed next to estimate times $t_1$ & $t_2$ and total recovery time from measurements. In the rest of the paper, a bold symbol signifies a complex quantity and subscript 'i' denotes the quantity at bus-i. The derivations in sub-section A & B are presented in more detail in [16] and are reproduced here for completeness.

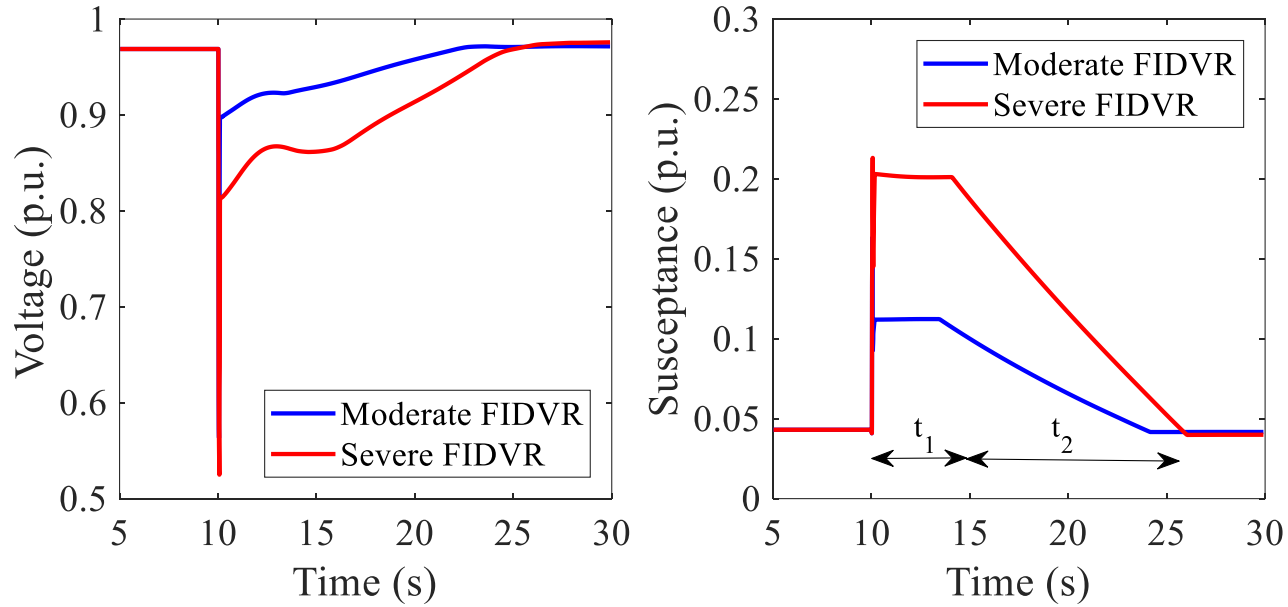

Fig. 4. Voltage response with various motor stalling proportion (left). Load susceptance with various motor stalling proportion (right). $t_1$ & $t_2$ are indicated for the severe FIDVR event.

### A. Time to initiate motor disconnection ($t_1$)

The expression for the load voltage $\boldsymbol{V_L}$ is given by (1). After the 1ϕ IM motor stalling, the thermal power dissipated is given by (2) and the corresponding differential equation for the temperature is given by (3). Initially the internal temperature is zero and the thermal loss is zero. As the stalling condition occurs suddenly, the input to the thermal delay block is a step function with value $P_{th}$ and the internal temperature increases exponentially as shown in (4). $G_{stall}$ remains same as $f_{th}$ is equal to 1 till the temperature reaches $\theta_1$ and the time taken for the temperature to reach $\theta_1$ can be calculated by substituting (2) in (4) to get (5).

$$\boldsymbol{V}_L = \boldsymbol{V}_{\mu PMU} - \boldsymbol{I}_L \cdot (R_f + j \cdot X_f) \quad (1)$$

$$P_{th} = V_L^2 \cdot G_{stall} \quad (2)$$

$$\frac{d\theta}{dt} = \frac{1}{T_{Th}}(P_{th} - \theta) \Rightarrow \theta = P_{th}\left(1 - e^{(-t/T_{th})}\right) \quad (3)$$

$$\theta_1 = P_{th}\left(1 - e^{(-t_1/T_{th})}\right) \quad (4)$$

$$t_1 \approx -T_{th} \cdot \ln\left(1 - \theta_1/(V_L^2 \cdot G_{stall})\right) \quad (5)$$

Next, we can determine the time taken for the motor temperature to rise from $\theta_1$ to $\theta_2$ by understanding how the thermal trip fraction $f_{th}$ varies with time.

### *B. Time to complete motor disconnection ($t_2$)*

The thermal trip fraction $f_{th}$ is a linear function of the internal temperature and is given by (6). To derive an expression describing how the thermal trip fraction varies with time, (6) is differentiated and the expression in (3) is substituted, leading to expression (7). In this expression, the voltage is an implicit function of the fraction $f_{th}$.

$$f_{th} = 1 - \frac{(\theta - \theta_1)}{(\theta_2 - \theta_1)} \quad (6)$$

$$\frac{d\, f_{th}}{dt} = \frac{-d\theta/dt}{(\theta_2 - \theta_1)} = \frac{\left(\theta_2 - (\theta_2 - \theta_1) \cdot f_{th} - V_L^2 \cdot G_{stall}\right)}{T_{th}(\theta_2 - \theta_1)} \quad (7)$$

Equation (7) enables us to understand the behavior of the thermal relay. Initially, the value of the voltage is low and the value of $f_{th}$ is 1, leading to a negative value of $df_{th}/dt$ and implying that the $f_{th}$ will reduce from 1 increasing the voltage magnitude. As the voltage increases and the $f_{th}$ decreases, the slope becomes further negative and increases the rate of rise of voltage. Finally, as the value of $f_{th}$ reaches 0, the voltage is close to the pre-contingency voltage at which time all the 1ϕ IM are disconnected and the thermal trip relay operation ends. The differential equation (7) is non-linear and can be numerically solved for a particular scenario but is difficult to analyze for a general case. In [16], we analyzed the differential equation (7) and derived an approximate expression for $t_2$ and the final expression is presented in (8) which needs the recovery voltage level at the bus when the FIDVR is complete. This is usually between 0.95 p.u. to 1 p.u.

$$t_2 \approx \frac{2T_{th}(\theta_2 - \theta_1)}{\left((V_L^2 + V_{L_{recover}}^2)G_{stall} - \theta_1 - \theta_2\right)} \quad (8)$$

$$t_{total} = t_1 + t_2 \quad (9)$$

The total time to recovery ($t_{total}$) is the sum of $t_1$ & $t_2$ and thus, the time to recover from FIDVR can be determined using the voltage and admittance measurements from the $\mu$PMU along with the load parameters. The recovery time is used an indicator for FIDVR severity. These expressions are used in the next section to formulate FIDVR mitigation using DER control and A/C control using smart thermostats.

## III. FIDVR Mitigation using DER and A/C Control

The FIDVR phenomenon is of concern to the system as sustained low voltages are not expected by the various components in the system and they might disconnect, leading to uncontrolled loss of generation/load. Usually, utilities need to satisfy a voltage recovery criteria (e.g. 0.95 p.u. in 10s [17]) that ensures that the system can recover after a fault. During FIDVR, this criteria can be violated (e.g. 0.95 p.u. in 13s) and so we need to take control to improve recovery time (by 3s in this case) to ensure the voltage satisfies the criteria. A change in the values $t_1$ & $t_2$ at bus-i can only occur due to a change in the voltage at bus-i and the relation between them is shown in (10). This expression is derived by linearizing the expressions in (5) and (8). The expressions for $dt_{1_i}/dV_{L_i}$ and $dt_{2_i}/dV_{L_i}$ can be analytically derived from (5) and (8).

$$\Delta t_{1_i} \approx \frac{dt_{1_i}}{dV_{L_i}} \cdot \Delta V_{L_i}; \Delta t_{2_i} \approx \frac{dt_{2_i}}{dV_{L_i}} \cdot \Delta V_{L_i} \quad (10)$$

Various control schemes at various locations in the network can lead to a change in the voltage at bus-i. Using linearization of the network equations, $\Delta V_{L_i}$ can be written as a linear combination of all possible controls scaled by the partial derivative as shown in (11). The total change in the recovery time due to the control throughout the network is then given by (12). The quantity $\partial V_{L_i}/\partial u_j$ is the change in the voltage at bus-i due to the m[th] control action and can be determined from information of the network parameters and the node voltages. These sensitivities need to be calculated during the early phase of FIDVR using the measurements from the $\mu$PMUs.

$$\Delta V_{L_i} \approx \sum_{j=1}^{m} \frac{\partial V_{L_i}}{\partial u_j} \Delta u_j \quad (11)$$

$$\Delta t_{total_i} = \Delta t_{1_i} + \Delta t_{2_i} \approx \left(\frac{dt_{1_i}}{dV_{L_i}} + \frac{dt_{2_i}}{dV_{L_i}}\right) \cdot \sum_{j=1}^{m} \frac{\partial V_{L_i}}{\partial u_j} \Delta u_j \quad (12)$$

Equation (12) can be written in a matrix form (13), where $D_{t_1,V}$ & $D_{t_2,V}$ are diagonal matrices of size $n \times n$ and the element $(i,i)$ is given by $dt_{1_i}/dV_{L_i}$ and $dt_{2_i}/dV_{L_i}$ respectively. $S_{V_L,u}$ is a matrix of sensitivities of size $n \times m$ and the element $(i,j)$ is given by $\partial V_{L_i}/\partial u_j$ and $\Delta u$ is a column vector of size $m \times 1$ which correspond to the various control schemes possible. The increasing number of active components in the DS provide the means to mitigate FIDVR. In this paper, we concentrate on the reactive power injection from PV inverters and A/C on/off control via smart thermostats [18]. Thus, the control vector $[\Delta u]$ can be written as $[\Delta u_{PV} \quad \Delta u_{AC}]^T$ and the expression for the change is recovery time is written as (14).

$$[\Delta t_{total}] = \left[D_{t_1,V} + D_{t_2,V}\right] \cdot \left[S_{V_L,u}\right] \cdot [\Delta u] \quad (13)$$

$$\begin{aligned}[\Delta t_{total}] &= \left[D_{t_1,V} + D_{t_2,V}\right] \cdot \left[S_{V_L,u_{PV}} \quad S_{V_L,u_{AC}}\right] \cdot \begin{bmatrix}\Delta u_{PV} \\ \Delta u_{AC}\end{bmatrix} \\ &= [A] \cdot \begin{bmatrix}\Delta u_{PV} \\ \Delta u_{AC}\end{bmatrix}\end{aligned} \quad (14)$$

All the elements in 'A' and $\Delta t_{total}$ are negative as Q-injection and load disconnection will reduce the recovery time. In this manner, we have derived a linearized expression to estimate the change in recovery time and this expression can be used to estimate the minimum control necessary to ensure voltage recovery within a specified amount of time. It is important to remember that the true system is non-linear and so this linearization is bound to have errors and this will be discussed in the results.

The control of PV devices essentially amounts to reactive power injection (assuming no curtailment of active power) and so the $S_{V_L,u_{PV}}$ matrix is estimated from the change in voltages due to reactive power injection at nodes with PV penetration. This is efficiently implemented as solving a set of linear equations from the power flow jacobian during the FIDVR event and is parallelized for multiple buses, speeding up the estimation of $S_{V_L,u_{PV}}$. Similarly, the control of A/C devices is equivalent to reducing the amount of active and reactive power demand at various nodes and so the $S_{V_L,u_{AC}}$ matrix is estimated from the change in voltages due to reducing active and reactive power at various nodes and can also be parallelized, ensuring that the $S_{V_L,u_{AC}}$ estimation is done in an online manner. As the matrices $D_{t_1,V}$ & $D_{t_2,V}$ have analytical expressions, their values are calculated very quickly and so the full matrix 'A' relating the control inputs to the change in recovery time can be calculated in an online manner and this is used to determine the effective regions and types of control for a specific FIDVR occurrence.

To determine the locations and amount of control, a linear optimization can be formulated using the linear relation (14). The formulation is shown in (15).

$$\begin{aligned} &\min c^T \cdot \Delta u \\ &-A \cdot \Delta u \geq -\Delta t_{total} \\ &\Delta u_{min} \leq \Delta u \leq \Delta u_{max} \end{aligned} \tag{15}$$

This formulation minimizes the control cost with coefficients 'c' and ensures that the recovery time improves by a minimum of $\Delta t_{total}$ while keeping the control within the bounds. The negative sign in the inequality constraint is present as all the elements in 'A' and $\Delta t_{total}$ are negative. $\Delta t_{total}$ is the change in the recovery time to ensure the recovery satisfies the voltage criteria and the constraints in (15) ensures that the voltage recovery improves at all the nodes.

One of the challenges in this approach is that $\mu$PMU measurements at all the nodes in the system are used to estimate 'A' in (15). This is not practical and so a methodology that requires lesser number of measurements is necessary. We propose to use the radial nature of the DS to aggregate the dynamic loads in an area into a reduced distribution system model that can be monitored and analyzed using lesser number of $\mu$PMU measurements. This has the added advantage solving the optimization problem (15) faster as the number of nodes are reduced from the original system. The methodology of the reduced distribution system model is described in detail in the next section.

## IV. Reduced Distribution System Model (RDSM)

In order to implement the FIDVR mitigation scheme that has been described in the previous section using few $\mu$PMUs, it is essential to aggregate a group of loads into a single load model. Consider a radial DS shown in Fig. 5 with N nodes with each node comprising of static, electronic, motor loads (3ϕ and 1ϕ AC motors) and PV inverters. Measurement devices such as μPMUs measure load voltage and power in the distribution lines/loads at sub-second intervals. The objective is to reduce the number of nodes and represent the load at each measurement node using an aggregated dynamic model that captures the overall dynamic behavior of the full model. The placement of μPMUs is a problem that is beyond the scope of this paper. For this paper, we assume that they are placed at nodes where secondary feeders and large loads are connected to the primary feeder. The proposed RDSM is made up of several sub-models connected in a structure similar to the original topology as shown in Fig. 6 (a).

The sub-model is analogous to the CLM described previously with selected parameters to represent relevant portions of the DS with an equivalent feeder impedance, a load tap changing transformer and a load block as shown in Fig. 6 (b). The load block includes static load, IM loads, and PV inverter. The static load parameters correspond to the conventional ZIP model. The 3ϕ IM (A, B, C type motors of the CLM) are lumped into one motor and the 1ϕ IM (Type D of the CLM) represents the motors used in residential A/C compressors.

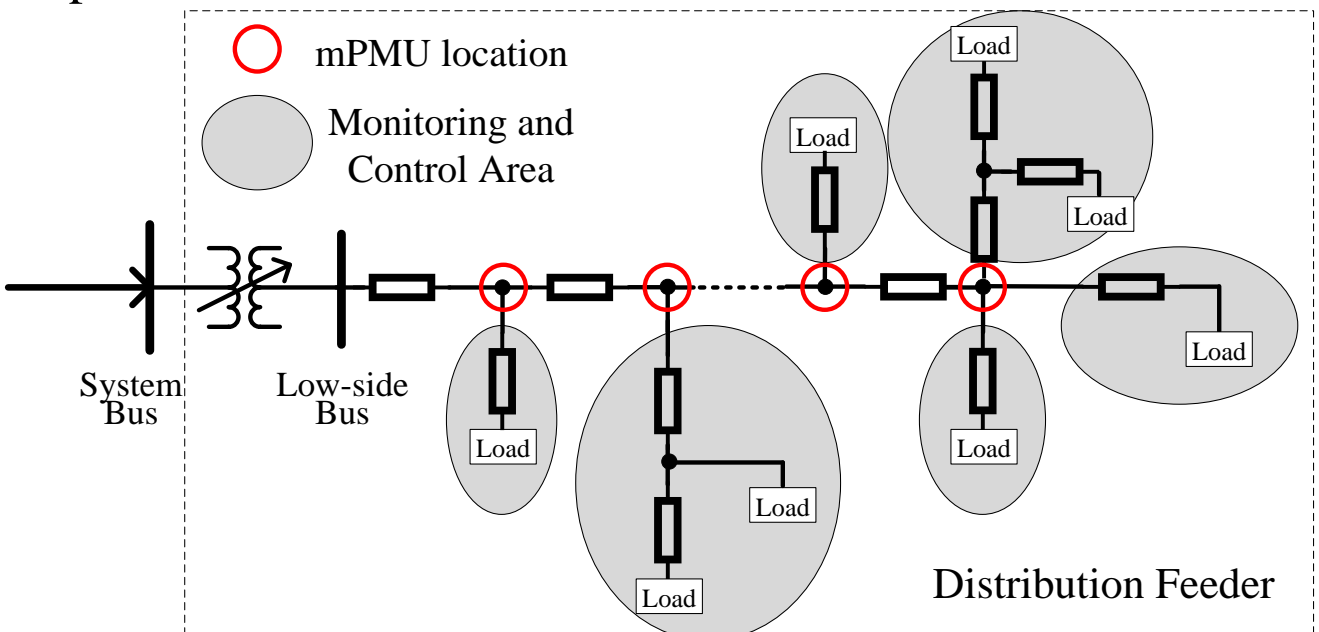

Fig. 5. Radial distribution system with μPMUS installed in some nodes

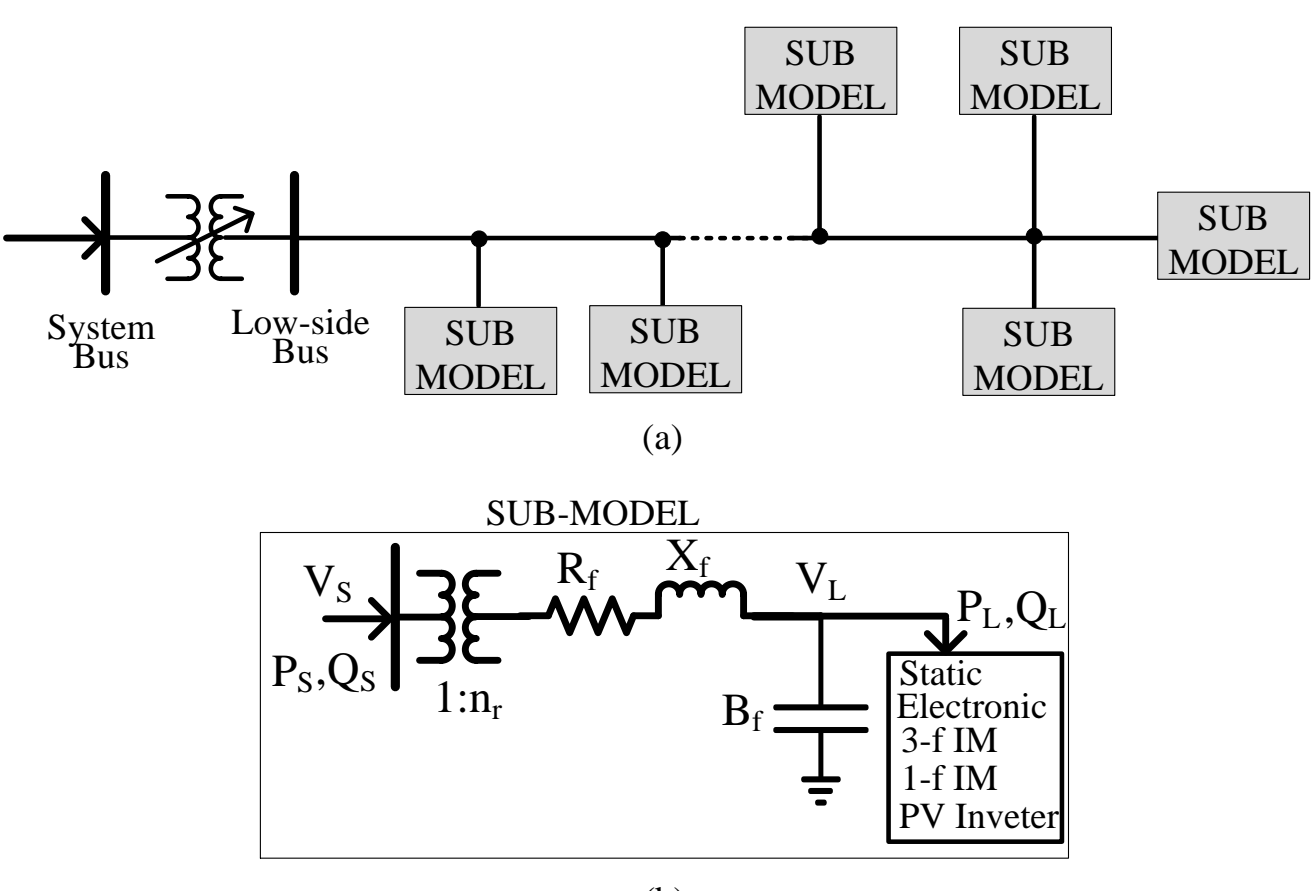

Fig. 6. (a) Generic reduced distribution system model. (b) Sub-model.

Table I shows the relevant parameters of the sub-model that would represent the portion of the DS network. Here, Fs, Fm1, Fm3 are the fractions of the corresponding loads and Fpv is the fraction of the equivalent PV in that portion of the network. Rf, Xf, Bf and nr are the parameters of the equivalent feeder impedance. The static load, the 3ϕ IM and the 1ϕ IM and PV inverter are represented by the parameters in the respective columns in Table I and these parameters are defined in the WECC CLM specifications [14].

Table I. Various load parameters of the RDSM

| Load Fraction | Equiv. Feeder | Static (ZIP) | 3ϕ IM | 1ϕ IM | PV |
|---|---|---|---|---|---|
| Fs | Rf | Pz0 | Rs | Vstall | Ppv |
| Fm1 | Xf | Qz0 | Xls | Tstall | Qpv |
| Fm3 | Bf | Pi0 | Xm | Rstall | |
| Fpv | nr | Qi0 | Rr1 | Xstall | |
| | | Pp0 | Xlr1 | $T_{th}$ | |
| | | Qp0 | Rr2 | $\theta_1$ | |
| | | Qsh0 | Xlr2 | $\theta_2$ | |
| | | | H | | |

### *A. RDSM Parameter Estimation using CoTDS simulation*

The Combined Transmission-Distribution System (CoTDS) simulation model [19] provides a means to generate surrogate data in the absence of a wide variety of data under various scenarios for the purpose of determining the parameters the DS. In [20], the CoTDS modeling methodology was used to calculate and validate the equivalent feeder impedance of the reduced order model of the of the entire DS using steady state data. This methodology is extended in this paper to determine the DS model parameters by also including the dynamics.

The dynamic data that is required is obtained by performing CoTDS simulation on a system comprising of a single-generator connected to the DS under study. The TS can either be a test system or an equivalent of a large TS under study. Since the purpose here is to generate large amount of surrogate data from the DS, it is not necessary to consider the entire TS. The entire system becomes necessary at a later stage when studying or validating the FIDVR control and mitigation scheme.

Fig. 7 shows the CoTDS Simulation set up for generating the measured data from the DS. The dynamic data is generated by applying a fault on the transmission-side. The fault is applied at different impedances and different fault times. This leads to various scenarios of TS fault which gives a variety of fault voltage and time for which the fault voltage is present. During the CoTDS simulation, the dynamics of each motor in the DS is uniquely calculated and their stalling condition is evaluated. As the fault scenario is varied, the number of motors stalling and their recovery is different and this leads to several data sets. The time series data for voltages, active power and reactive power is recorded at the nodes where the μPMUs are placed. The data thus generated is used for determining the parameters of the sub-models of the DS.

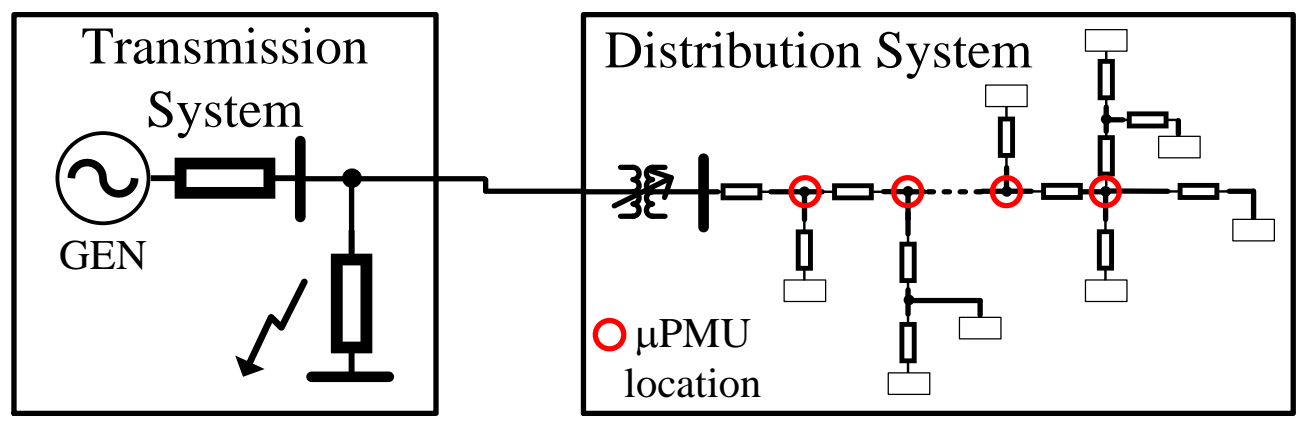

Fig. 7. CoTDS Simulation set up to generate the dynamic measurement data

### *B. Determination of RDSM parameters using non-linear optimization*

The sub-model parameters are classified into steady-state network parameters and dynamic load parameters. The steady-state network parameters correspond to the equivalent feeder impedance and the dynamic parameters correspond to the load component parameters as given in Table I. In the absence of real measured data, surrogate data obtained from a CoTDS simulation is used. The steps to determine these parameters are:

1. The CoTDS simulation is run using a single generator and a single line TS and the DS that needs to be reduced.
2. The steady state data of the sub-station voltage, the active power, reactive power and the voltage data at all the nodes of the distribution feeder are used to determine the equivalent feeder parameters using the method described in [20].
3. The dynamic data of the voltage, the active power and the reactive power at the transmission side is recorded. This represents the actual data from the actual load.
4. The values of the parameter set $\lambda$ of the sub-model are determined using an optimization routine to minimize the error between the time series of the measured data, $D$, and the calculated values, $C(\lambda)$. Fig. 8 shows the processing of the data in an optimization routine to estimate the parameters of the sub-model. The objective function of the optimization that needs to be minimized, $\eta(\lambda)$, is the sum of squares of the difference between the two time series and is given by (16):

$$\eta(\lambda) = [D - C(\lambda)]^T . [D - C(\lambda)] \qquad (16)$$

where $C(\lambda)$ is the corresponding calculated values of the data set, ***D*****,** for a given $\lambda$. The calculated values $C(\lambda)$ are obtained by solving the dynamic equations of the sub-model including the effects of the stalling and thermal tripping of the 1ϕ IM.

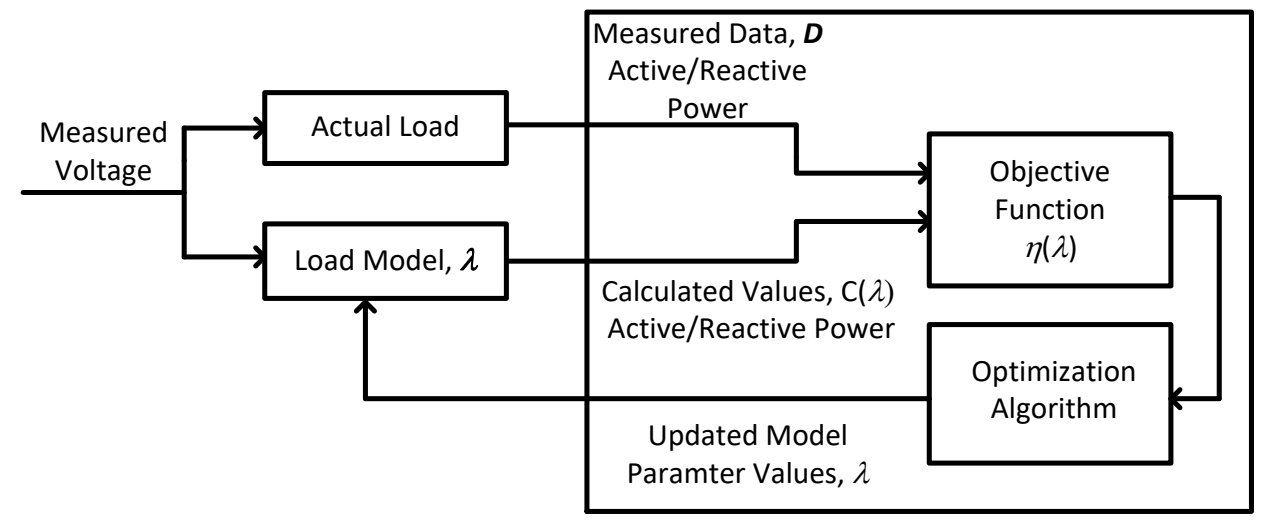

Fig. 8. Block diagram of methodology for RDSM parameters determination

The estimated RDSM parameters can then be used for FIDVR monitoring and mitigation and Fig. 9 summarizes the proposed methodology using RDSM and $\mu$PMU measurements. Next numerical results validating proposed methodology are discussed using an example system.

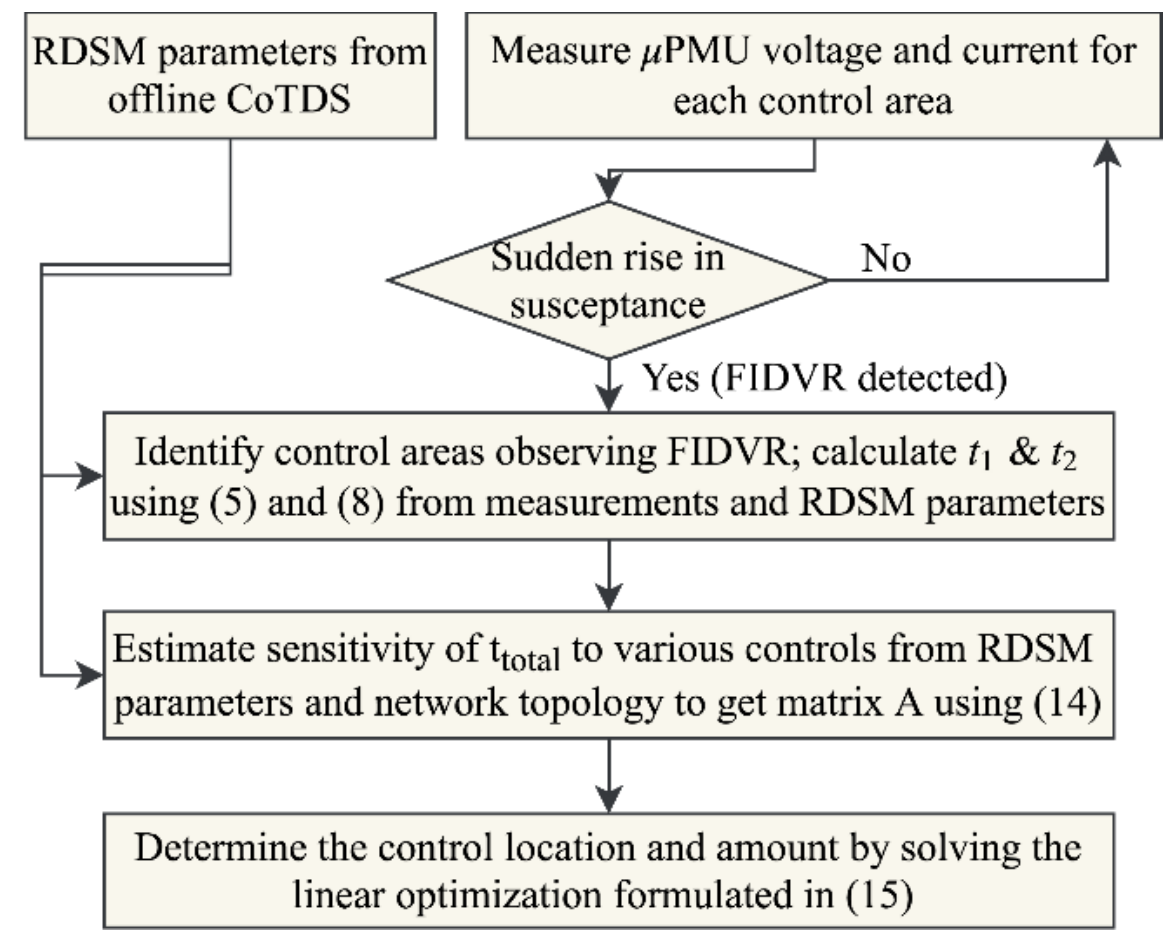

Fig. 9. Flowchart for detecting and monitoring FIDVR using measurements

## V. NUMERICAL RESULTS ON TEST SYSTEM

An IEEE 9-bus TS and an IEEE 37-node DS [21] is considered for implementing the control scheme that is developed to mitigate the FIDVR. Fig. 10 shows the interconnection of these systems with the IEEE 37-node system connected to Bus 6 of the IEEE 9-bus TS. The power demanded by the DS is 2.5 MW and so the remaining power at the transmission bus is assigned to other feeders that are not under study. The DS is shown in Fig. 10 & is divided into 6 load areas with the root node voltage and the currents in each area being measured by a μPMU (in red) as shown in Fig. 7. Observe that a few areas share the same root node and the $\mu$PMU located here should measure the currents into load area separately.

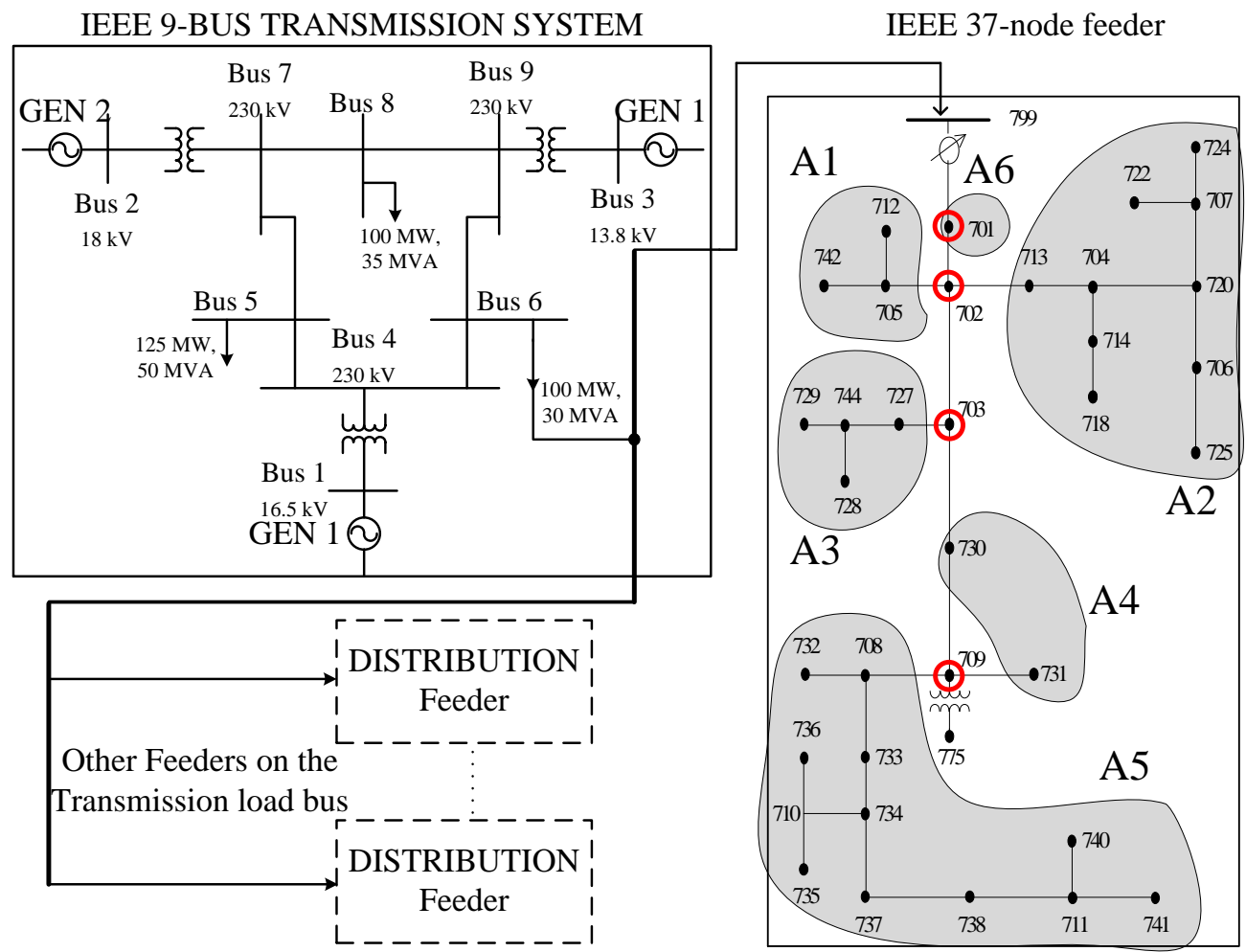

Fig. 10. Test system used to validate proposed methodology. The areas in the distribution system are shaded and the red nodes are the $\mu$PMU locations.

## A. RDSM parameters of IEEE 37-node distribution system

Each of the DS loads in the IEEE 37-node feeder is separated into the composite load model components including static, electronic, 3ϕ IM and 1ϕ A/C IM. In order to simulate a realistic scenario, the fraction of loads of each type (*Fs, Fel, Fm*3 *and Fm*1) is assigned according a normal distribution around a mean value which is estimated based on the type of loads (residential, industrial or commercial) present in each location [22]. In addition, each of the motor load types which have their own set of parameters to characterize them and have variability included by connecting several motors with a normal distribution of parameters. This procedure leads to a comprehensive and detailed model of the DS.

To test the system behavior, a fault is created at node 701 in the DS for a duration of 80 ms and the voltages observed by the $\mu$PMUs and the substation are plotted in Fig. 11. It can be seen that FIDVR is observed in all the $\mu$PMUs due to the high proportion of the 1ϕ IM in the feeder. As the FIDVR occurs on only a single distribution feeder whose load (2.5 MW) is small compared to the load at the TS (100MW), the TS is not impacted by this fault and this is reflected in the substation voltage being nearly flat during the FIDVR. These voltage profiles are similar to the FIDVR data from $\mu$PMUs in the Southern California Edison system [10].

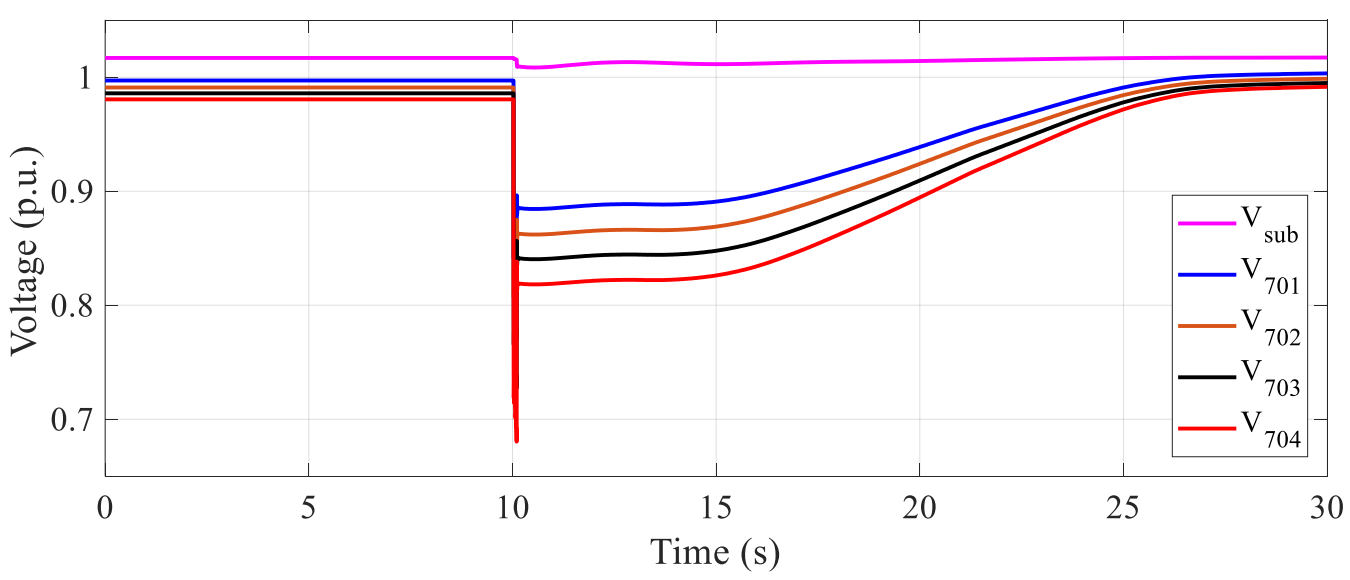

Fig. 11. Voltages at the substation and at $\mu$PMUs for a fault at node 701

The procedure described in section IV-B for estimating RDSM parameters is applied to the IEEE 37-node DS. From the different sets of data obtained from CoTDS, the sub-model parameters for each of the 6 areas are determined according to the optimization procedure described in the section IV-B. A few parameters for each of the control area are listed in Table II.

Table II. Sub-model parameters of the load areas

| Parameter | Root Node and Sub-Model Parameter Values | | | | | |
|---|---|---|---|---|---|---|
| | A1 | A2 | A3 | A4 | A5 | A6 |
| Root Node | 702 | 702 | 703 | 709 | 709 | 701 |
| $P_{load}$ (kW) | 178 | 538 | 245 | 160 | 684 | 420 |
| Fs | 0.61 | 0.46 | 0.49 | 0.49 | 0.47 | 0.2 |
| Fm1 | 0.39 | 0.54 | 0.29 | 0.51 | 0.53 | 0.1 |
| Fm3 | 0 | 0 | 0.22 | 0 | 0 | 0.7 |
| Rstall | 0.061 | 0.092 | 0.057 | 0.074 | 0.072 | 0.080 |
| Xstall | 0.073 | 0.112 | 0.058 | 0.077 | 0.091 | 0.090 |
| $T_{th}$ | 17.84 | 12.00 | 15.14 | 13.99 | 13.62 | 15.00 |
| $\theta_1$ | 0.714 | 0.452 | 0.450 | 0.653 | 0.739 | 0.800 |
| $\theta_2$ | 3.025 | 1.949 | 3.750 | 3.222 | 2.615 | 3.000 |

## B. Validation of RDSM parameters with CoTDS simulation

To demonstrate the accuracy of the RDSM model behavior compared to the full model, the active and reactive powers of the different areas are plotted in Fig. 12 and Fig. 13 for a fault scenario for different areas using both the CoTDS simulation and the RDSM. It can be seen that the active and reactive power profiles of the RDSM closely matches the data from CoTDS simulation at all the load areas for almost the entire recovery period after the fault and captures the FIDVR behavior. The active and reactive powers also matched well for various other fault scenarios. These plots and observations verify that the RDSM can indeed capture the overall behavior of the full model with a reasonable degree of accuracy and validates the parameters of the RDSM. Next, the accuracy of the proposed FIDVR monitoring and mitigation methodology is tested on this system. In the rest of the paper, the recovery time is defined as the time taken for the voltage at all the $\mu$PMU locations to recover to 0.95 p.u.

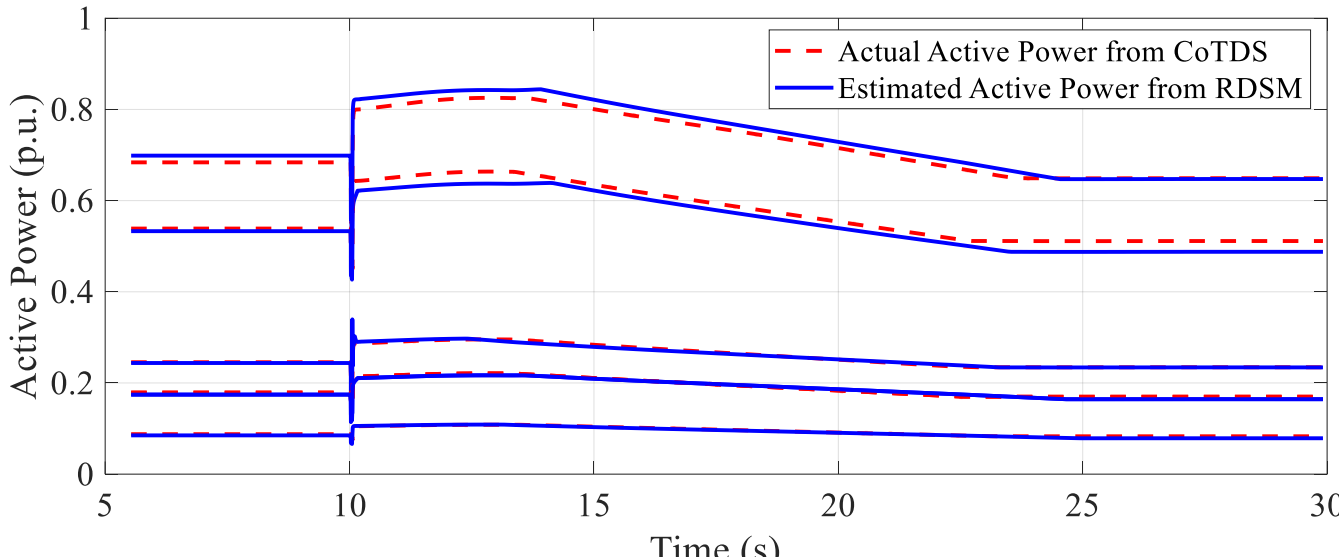

Fig. 12. Active power of the areas using CoTDS and the RDSM parameters

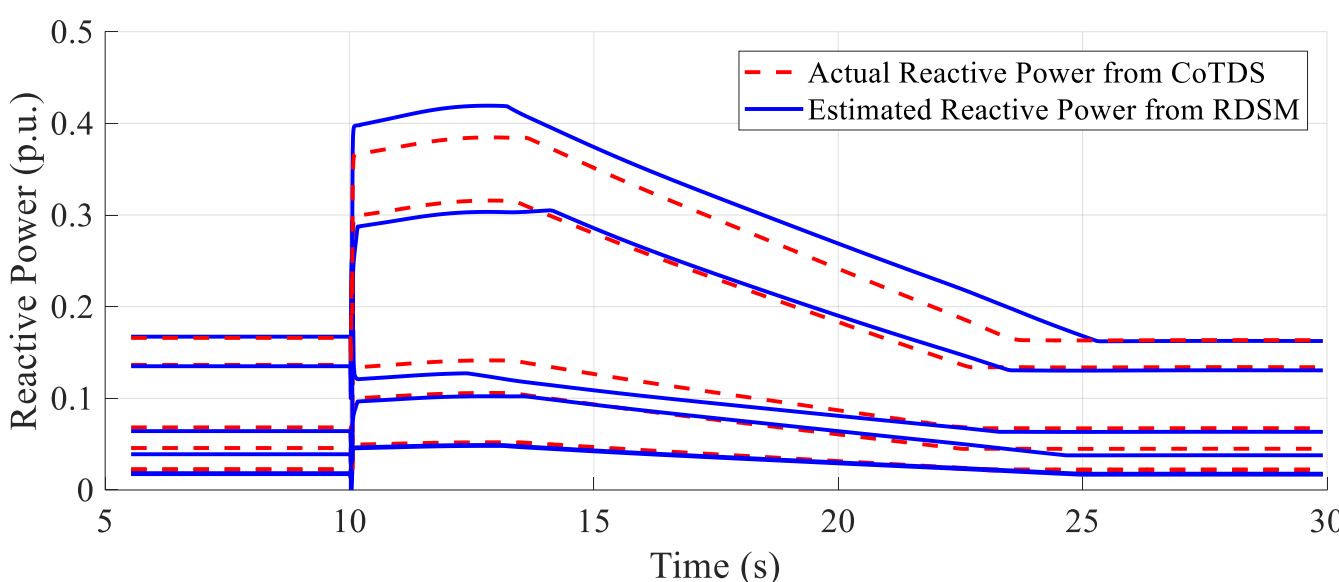

Fig. 13. Reactive power of the areas using CoTDS and the RDSM parameters

## C. Recovery time estimation for monitoring FIDVR

To validate the expressions in (5) and (8), various faults in the DS are created in the full CoTDS simulation and the actual recovery time is measured from the resulting FIDVR profile. This is compared to the estimated recovery time calculated using (9) and the $\mu$PMU measurements at the root nodes of areas A1 – A6. Fig. 14 plots the voltage response at node 709 in the DS for the three fault locations with varying fault duration applied in various areas of the distribution feeder. The recovery times are presented in Table III and they demonstrate that the estimated recovery time lies within 15 % of the actual time in all the cases with the largest errors occurring in scenarios with low fault durations. The recovery time can be estimated in <0.5s after the FIDVR event is detected and enables the fast detection of events that are likely to exceed the recovery time specified by the utility (e.g. 10s). The fast detection and recovery time estimation of FIDVR makes it possible to initiate control schemes to improve recovery time and this is described next.

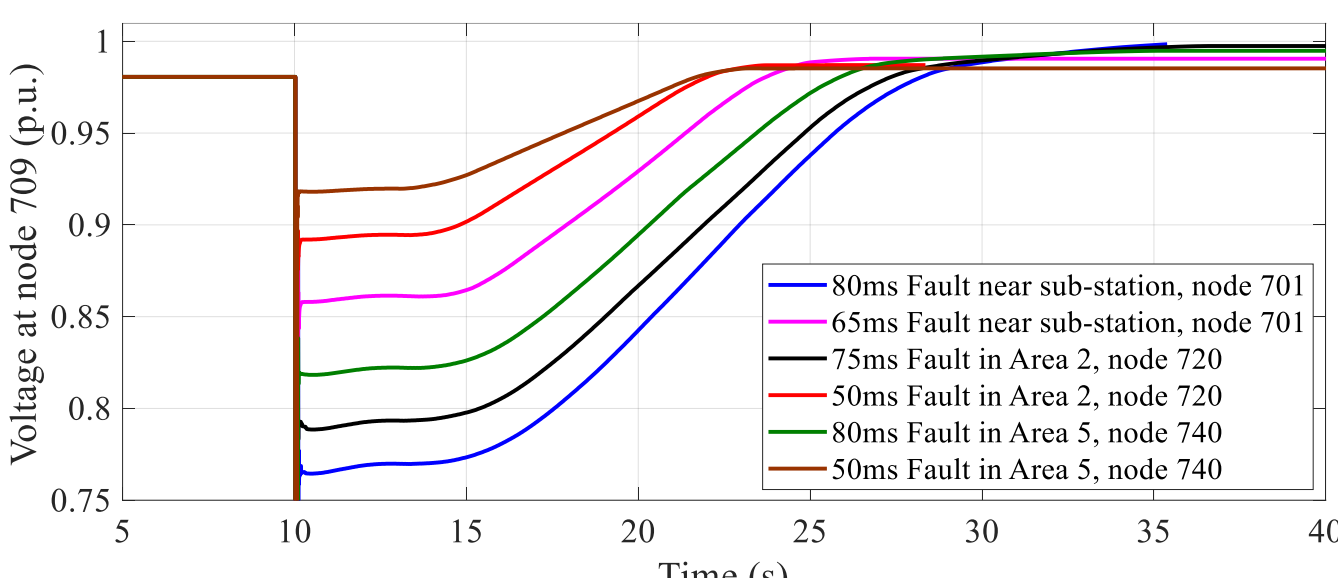

Fig. 14. Voltage at node 709 for various faults in the DS

Table III. Comparison between the actual and estimated recovery times

| Fault location | Fault duration | Actual $t_{total}$ | Estimated $t_{total}$ | Abs. Error (%) |
|---|---|---|---|---|
| 701 (near substation) | 80 ms | 15.7 s | 14.9 s | 5 % |
| | 65 ms | 11.4 s | 12.1 s | 6 % |
| 720 (in A2) | 75 ms | 14.7 s | 14.1 s | 4 % |
| | 50 ms | 9.2 s | 10.2 s | 11 % |
| 740 (in A5) | 80 ms | 13.6 s | 13.4 s | 1.5 % |
| | 50 ms | 7.9 s | 9.0 s | 13.5% |

### D. *Control for mitigating FIDVR*

To verify that the expression in (14) can predict the FIDVR recovery time improvement, the case with a fault at bus 740 for 80 ms is chosen. FIDVR is detected by the sudden rise in susceptance and the sensitivities are calculated at the FIDVR condition from measurements, topology & RDSM parameters. The controls (A/C load disconnection and reactive power injection from PV inverters ($f_{PV}$ = 25%)) are implemented in the full CoTDS and are triggered 1.5s after the FIDVR is detected. This time delay is to account for the communication delays and computation time to estimate sensitivities and execute the optimization. As there are no measurements within each load area to identify the particular motors that are stalled, a specified percentage of motors are randomly disconnected in each area. This is the practical scenario as we cannot identify the individual stalled motors. The recent IEEE 1547 [23] standard mandates that new PV inverters should be able to inject reactive power corresponding to 44% of its rating without active power curtailment and is implemented in the CoTDS for reactive power control. The various control scenarios are listed in Table IV along with the actual and estimated $\Delta t_{total}$. Fig. 15 plots the voltage at node 709 for the various control scenarios.

It can be seen from the results that the estimated $\Delta t_{total}$ for various controls match the actual $\Delta t_{total}$ from the CoTDS. This validates the derivation of the change in recovery time using sensitivities in (14). However, as the phenomenon is inherently non-linear and the sensitivities are a linearized representation, it is expected that as the control amount increases, the error between the actual and estimated values will increase and this is precisely what is observed from Table IV.

Table IV. Comparison between the actual and estimated recovery time improvement for various control actions for fault at node 740

| Control Description | Load disconnected | Actual $\Delta t_{total}$ | Estimated $\Delta t_{total}$ |
|---|---|---|---|
| 10% A/C disconnection in all areas; no Q from PV | 91.7 kW | -0.95 s | -0.9 s |
| 20% A/C disconnection in all areas; no Q from PV | 183.4 kW | -2.05 s | -1.85 s |
| 30% A/C disconnection in all areas; no Q from PV | 275.1kW | -3.40 s | -2.8 s |
| 30% A/C disconnection in Area 5; no Q from PV | 108.8 kW | -2.10 s | -1.8 s |
| 44% Q-Injection from PV ($f_{PV}$=0.25) in all areas | - | -0.65 s | -0.7 s |

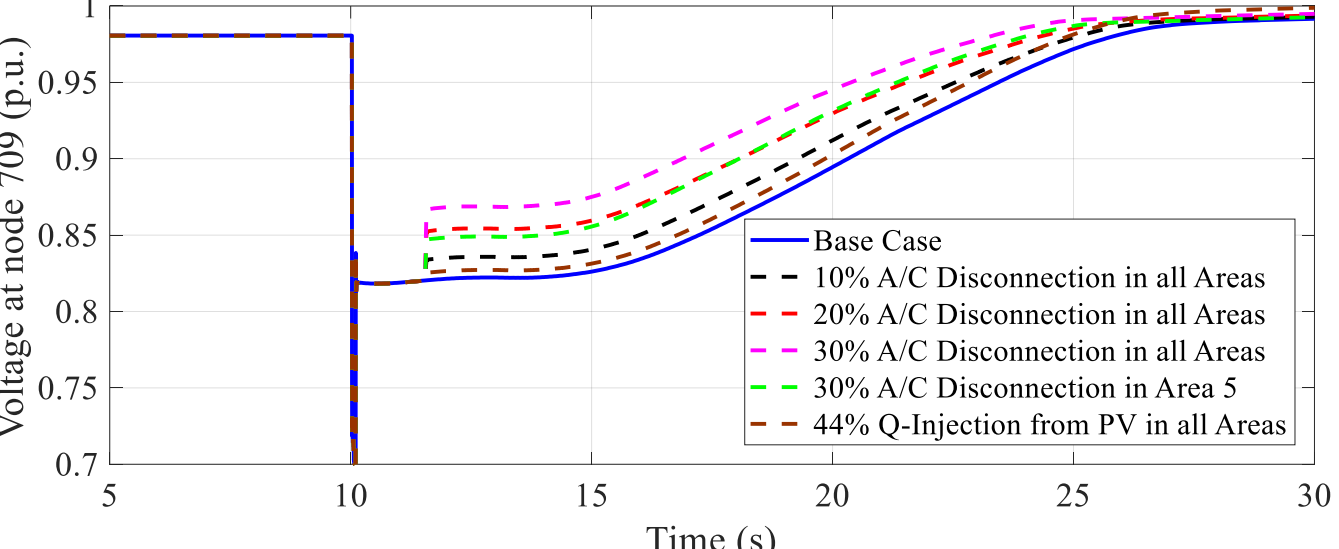

Fig. 15. Voltage at node 709 for various control actions for fault at node 740

It can also be seen that the same improvement in the recovery time can be realized by disconnecting 108.8 kW of A/C load in area 5 versus 183.4 kW of A/C load in the overall DS – A reduction of 40% in the load disconnected. This is also captured by the estimated $\Delta t_{total}$ as the calculated sensitivities of the control in area 5 are the highest in magnitude. Thus, the sensitivity based method can estimate the $\Delta t_{total}$ and can be used for determining effective control locations. Next, the optimization formulation (15) is implemented in Matlab to estimate the minimum A/C disconnection to improve voltage recovery by various $\Delta t_{total}$ values. A control constraint of 50% disconnection in each area is imposed for A/C load control. The estimation of the sensitivities and the execution of the linear optimization took <0.5s. Two control limits for the PV inverters were evaluated – normal unity power factor and 44% maximum Q-injection. The resulting optimal control schemes for a 3 scenarios with varying $\Delta t_{total}$ are listed in Table V and the voltage responses at node 709 are plotted in Fig. 16.

For the 1st case, the $\Delta t_{total}$ is short enough so that control in A5, which has the highest sensitivity, is sufficient to satisfy the optimization constraints. For the 2nd case, the control limit in the A5 is reached and the optimization selects the next sensitive load area, A2, for control. The $\Delta t_{total}$ in this case corresponds to 30% load disconnection in all areas in Table IV. It can be observed that the amount of load disconnection dropped from 275 kW to 200 kW - a reduction of 28%, demonstrating the benefit of the proposed methodology. If the unity power factor constraint is relaxed and the PV inverters are allowed to inject reactive with no active power curtailment, the load disconnection is further reduced to 145 kW – an overall reduction by 47%, validating the utility of controlling reactive injection from DERs. While we have concentrated on one fault scenario here, the proposed methodology is able quickly (<0.5s) calculate the A/C load disconnection and Q-injection from DERs to mitigate FIDVR due to faults in various locations with similar reduction in A/C load disconnection compared wide spread control in the overall distribution system.

Table V. Comparison between various control actions to improve the voltage recovery with different control constraints for fault at node 740

| $u_{max}$ Constraints | $\Delta t_{total}$ | Optimal Control Description | Load disconn. | Actual $\Delta t_{total}$ |
|---|---|---|---|---|
| No PV-Q<br>50% A/C | -2 s | A5 - 33% A/C load | 120 kW | -2.35 s |
| No PV-Q<br>50% A/C | -3.4 s | A3 - 30% A/C load<br>A5 - 50% A/C load | 200 kW | -3.65 s |
| 44% PV-Q<br>50% A/C | -3.4 s | Full Q from PV<br>A5 - 30% A/C load | 145 kW | -3.5 s |

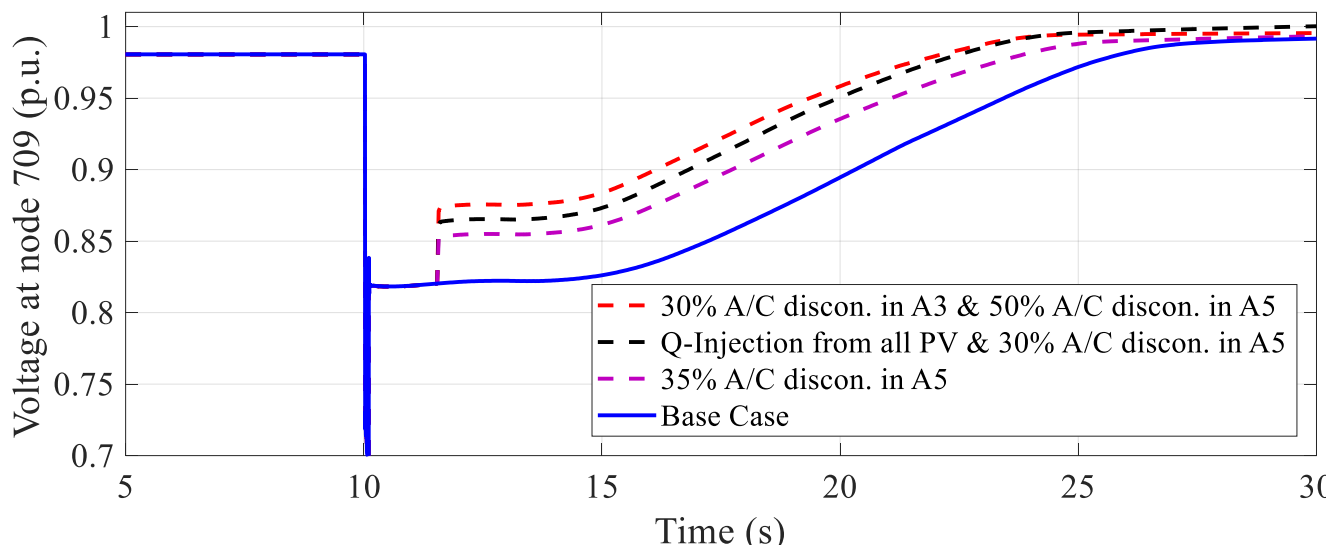

Fig. 16. Impact of optimal control to improve voltage for fault at node 740

With these results, we have verified the accuracy of the proposed methodology and demonstrated the reduction in load disconnection by utilizing the proposed control scheme and thus, the paper can be concluded.

## VI. Conclusion

In this paper, a methodology to monitor and mitigate the phenomenon of FIDVR in the DS is proposed using PMU measurements. The expressions for recovery time are derived by simplifying the composite load model during FIDVR. A linearization of the non-linear expressions is then done to estimate the change in recovery time under various control schemes and a linear optimization problem is formulated to estimate the minimal control action necessary to recover faster.

In order to apply this method to a DS with few μPMUs, the Reduced Distribution System Model is proposed which is composed of sub-models that are analogous to the WECC CLM and aggregates the DS into load areas while ensuring the overall dynamics are retained. To test the proposed scheme, a dynamic co-simulation is performed with several fault scenarios on the IEEE 37 node DS connected to IEEE 9 bus TS. This RDSM is shown to capture the dynamic behavior of the full distribution system under various fault scenarios. The optimal control actions calculated by using the linear sensitivities are quickly (<0.5s) able to identify the critical regions for control and demonstrate that optimal control reduces the amount of load control significantly (>25%). Furthermore, Q-injection by DERs can be incorporated into the optimization, further reducing the load control (>40%). Thus, the proposed methodology enables online monitoring and mitigation of FIDVR by utilizing Q-injection from DERs with minimal load disconnection making it a promising application of $\mu$PMU measurements to enhance operation of the distribution systems.

Incorporating the non-linearity of FIDVR behavior into the optimal control estimation will reduce the error due to linearization and is the next step in our research. Also, analyzing the circumstances under which curtailment of active power and increasing Q-injection improves the voltage profile is a very relevant research direction. This will enable us to exploit the capabilities of DERs to improve distribution grid behavior during emergencies.